\documentclass{amsart}

\usepackage{amssymb}

\usepackage{graphicx}

\newtheorem{theorem}{Theorem}[section]
\newtheorem{lemma}[theorem]{Lemma}

\theoremstyle{definition}
\newtheorem{definition}[theorem]{Definition}
\newtheorem{example}[theorem]{Example}

\newtheorem{proposition}[theorem]{Proposition}
\newtheorem{corollary}[theorem]{Corollary}

\newtheorem{question}[theorem]{Question}

\theoremstyle{remark}
\newtheorem{remark}[theorem]{Remark}

\newcommand{\bC}{\ensuremath{\mathbb{C}}}

\newcommand{\bK}{\ensuremath{\mathbb{K}}}

\newcommand{\bQ}{\ensuremath{\mathbb{Q}}}
\newcommand{\bR}{\ensuremath{\mathbb{R}}}

\newcommand{\bZ}{\ensuremath{\mathbb{Z}}}

\newcommand{\scA}{\ensuremath{\mathcal{A}}}
\newcommand{\scB}{\ensuremath{\mathcal{B}}}

\newcommand{\scE}{\ensuremath{\mathcal{E}}}

\numberwithin{equation}{section}

\newcommand{\Der}{\operatorname{Der}}

\begin{document}

\title[Multiarrangements]{On the extendability of free multiarrangements}



\author{Masahiko Yoshinaga}
\address{Department of Mathematics, Faculty of Science, 
Kobe University, Kobe 657-8501, Japan}
\curraddr{}
\email{myoshina@math.kobe-u.ac.jp}

\thanks{}

\subjclass[]{}

\date{\today}

\dedicatory{}

\begin{abstract}
A free multiarrangement of rank $k$ is defined to 
be extendable if it is obtained from a simple rank $(k+1)$ 
free arrangement by the natural restriction 
to a hyperplane (in the sense of Ziegler). 
Not all free multiarrangements are extendable. 
We will discuss extendability of 
free multiarrangements for a special class. 
We also give two applications. 
The first is to produce 
totally non-free arrangements. 
The second is to give interpolating 
free arrangements between extended Shi 
and Catalan arrangements. 
\end{abstract}

\maketitle


\section{Introduction}

Let $V=\bC^\ell$ be a complex vector space with 
coordinate $(x_1, \cdots, x_\ell)$, $\scA=\{H_1, \ldots, H_n\}$ 
be a central arrangement of hyperplanes. Let us denote by 
$S=\bC[x_1, \ldots,x_\ell]$ the polynomial ring and fix 
$\alpha_i\in V^*$ a defining equation of $H_i$, i.e., 
$H_i=\alpha_i^{-1}(0)$. 
A multiarrangement is a pair $(\scA, m)$ of an arrangement 
$\scA$ with a map $m:\scA\rightarrow \bZ_{\geq 0}$, called 
the multiplicity.  
We denote 
$Q(\scA, m)=\prod_{i=1}^n\alpha_i^{m(H_i)}$ and $|m|=\sum_i m(H_i)$. 
An arrangement $\scA$ can be identified with a multiarrangement 
with constant multiplicity $m\equiv 1$, which is sometimes called 
a simple arrangement. 
Under these notation, the main object in this article is 
the following module of $S$-derivations which has contact 
to each hyperplane in the order $m$. 
\begin{definition}
Let $(\scA, m)$ be a multiarrangement, then define 
$$
D(\scA, m)=\{\delta\in\Der_S| \delta\alpha_i\in(\alpha_i)^{m(H_i)},
\forall i\}. 
$$
\end{definition}
The module $D(\scA, m)$ is obviously a graded $S$-module. 
A multiarrangement $(\scA, m)$ is said to be free with 
exponents $(e_1, \ldots, e_\ell)$ if 
$D(\scA, m)$ is an $S$-free module and there exists a 
basis $\delta_1, \ldots,\delta_\ell\in D(\scA, m)$ such that 
$\deg \delta_i=e_i$. Note that the degree 
$\deg\delta$ of a derivation $\delta$ is the polynomial 
degree, that is defined by $\deg (\delta f)=\deg\delta + 
\deg f-1$ for a homogeneous polynomial $f$. 
An arrangement $\scA$ is said to be free if $(\scA, 1)$ is free. 
Here we recall that the freeness is closed under localization. 
More precisely, let $X\subset V$ be a subset and define 
$\scA_X=\{H\in\scA \mid H\supset X\}$. Then the freeness of 
$(\scA, m)$ implies that of $(\scA_X, m|_{\scA_X})$.

A multiarrangement naturally appears as a 
restriction of a simple arrangement \cite{zie-multi}. 
Let $\scA$ be an arrangement. 
The arrangement $\scA$ determines the restricted 
arrangement $\scA^H=\{H\cap H'\mid H'\in\scA, H'\neq H\}$ 
on $H\in\scA$. The restricted arrangement $\scA^H$ possesses a 
natural multiplicity 
$$
\begin{array}{cccl}
m^H:&\scA^H&\longrightarrow&\bZ\\
&X&\longmapsto&\sharp\{H'\in\scA\mid X=H\cap H'\}. 
\end{array}
$$
The freeness of $\scA$ and $(\scA^H, m^H)$ 
are connected by the following theorem due to Ziegler. 

\begin{theorem}
\label{thm:ziegler}
\cite{zie-multi}
If $\scA$ is free with exponents $(1, e_2, \ldots, e_\ell)$, 
then the restriction $(\scA^H, m^H)$ is free with 
exponents $(e_2, \ldots, e_\ell)$. 
\end{theorem}
Recently freeness of multiarrangements are extensively 
studied \cite{atw-char, atw-emulti, st-double, ter-multi, waka-exp, wake-yuz}. 
The motivation to this article is to ask whether if a free multiarrangement 
is obtained as a restriction of a free simple arrangement. 
Theorem \ref{thm:ziegler} leads us to introduce the following notion, 
which seems to 
give an important class of free multiarrangements. 
\begin{definition}
Let $(\scA, m)$ be a free multiarrangement in $\bK^\ell$. 
We say $(\scA, m)$ is extendable if it can be obtained 
as a restriction of a free simple arrangement in $\bK^{\ell+1}$. 
\end{definition}

\begin{example}
(Non-extendable free multiarrangement) 
Consider a multiarrangement in $\bR^2$ 
$$
Q(\scA, m)=x^3 y^3 (x-y)^1 (x-\alpha y)^1 (x-\beta y)^1, 
$$
with $\alpha, \beta \neq 0,\pm 1$ and assume $\alpha$ and $\beta$ are 
algebraically independent over $\bQ$. (Indeed $\alpha\beta\neq 1$ 
is enough.)  
If the slopes $\alpha$ and $\beta$ are generic, 
then $(\scA, m)$ is 
free with exponents $(4,5)$ \cite{wake-yuz}. 
So the product of exponents is always $\leq 20$. 
We can prove that 
it is not extendable. It can be proved as follows 
(details are left to the reader). 
The deconing $\overline{\scA}$ (\cite{ot}) with 
respect to the hyperplane at infinity 
of an extension of $(\scA, m)$ is an affine 
line arrangement $\bR^2$ having the following defining equations: 
\begin{eqnarray*}
x&=&a_1, a_2, a_3, \\
y&=&b_1, b_2, b_3, \\
x-y&=&c, \\
x-\alpha y&=&d, \\
x-\beta y&=&e, 
\end{eqnarray*}
where $a_i, b_i, c, d, e\in\bR$. 
The characteristic polynomial $\chi(\overline{\scA}, t)$ is of the form 
$\chi(\overline{\scA}, t)=t^2-9t+p$, and we can prove that $p>20$. Thus 
$\chi(\overline{\scA}, t)$ is not factored. It follows from 
Terao's factorization theorem (\cite{ter-fact}) that 
any extension of $\scA$ is not free. 
\end{example}
Thus a free multiarrangement $(\scA, m)$ is not 
necessarily extendable in general. In the next section, 
we focus on some special kind of multiarrangements. 

\section{Extendability of locally $A_2$ arrangements}

\begin{definition}
An arrangement $\scA=\{H_1, \ldots, H_n\}$ 
is said to be {\em locally $A_2$} if $|\scA_X|\leq 3$ is satisfied 
for any codimension two intersection $X$. 
A system of defining equations $\{\alpha_1, \ldots, \alpha_n\}$ 
of a locally $A_2$ arrangement $\scA$ is called a 
{\em positive system} if the following condition 
is satisfied: 
Suppose $X$ is a codimension two intersection 
with $|\scA_X|=3$. Setting $\scA_X=\{H_i, H_j, H_k\}$. 
Then one of 
$\alpha_i=\alpha_j+\alpha_k$, 
$\alpha_j=\alpha_i+\alpha_k$ or 
$\alpha_k=\alpha_i+\alpha_j$ holds. 
\end{definition}

\begin{example}
The following are examples of locally $A_2$ arrangements with positive 
systems. 
\begin{itemize}
\item [(1)] Generic in codimension three. (Equivalently, $|\scA_X|=2$ for 
any codimension two intersection $X$.) In this case any system of 
defining equations is a positive system. 
\item [(2)] Coxeter arrangement of type $ADE$. 
In this case, a positive root system is corresponding to a positive system 
of defining equations. 
\item [(3)] Subarrangements or direct products of locally $A_2$ 
arrangements with positive systems possess the same property. 
Especially, this class is closed under localization. 
\item [(4)] (Shi arrangement of type $A_2$) $Q=xyz(x+y)(x-z)(y-z)(x+y-z)$. 
\end{itemize}
\end{example}

\begin{remark}
Note that a locally $A_2$ arrangement does not necessarily 
have a positive system (e. g., $Q=xyz(x+y)(x-z)(y-z)(x+y-2z)$). 
\end{remark}

We will discuss the extendability for multiarrangements 
of this class. 
More precisely, we consider the following concrete extension 
$E(\scA, m)$ of $(\scA, m)$ for given 
locally $A_2$ arrangement $\scA$ with a positive 
system $(\alpha_H)_H$. 
Let $(x_1, \ldots, x_\ell, z)\in\bC^\ell\times\bC$ be 
a coordinate system of $V\times\bC$ and 
define 
$$E(\scA, m)=\{z=0\}\cup 
\left\{
\alpha_H=kz
\left| 
k\in\bZ, -\frac{m(H)-1}{2}\leq k\leq \frac{m(H)}{2}
\right.\right\}.
$$
Then, if we denote $H_0=\{z=0\}$, 
it is obvious that $(E(\scA, m)^{H_0}, m^{H_0})=(\scA, m)$. 
Let us define the deconing of $E(\scA, m)$ as follows: 
$$\mathbf{d}E(\scA, m)=
\left\{
\alpha_H=k
\left| 
k\in\bZ, -\frac{m(H)-1}{2}\leq k\leq \frac{m(H)}{2}
\right.\right\}.
$$
Note that $\mathbf{d}E(\scA, m)$ is an affine 
arrangement in $V$.  

\begin{remark}
The above definition is motivated by 
that of the extended Catalan and Shi arrangements 
\cite{ede-rei}. 
Indeed, let $\scA$ be a Coxeter arrangement 
of type ADE. Choose the positive root system as 
the positive system as above. 
For a given positive integer $k\in\bZ_{>0}$, 
consider constant multiplicities 
$m=2k$ and $m=2k+1$. Then $E(\scA, 2k+1)$ is so 
called the extended Catalan arrangement and 
$E(\scA, 2k)$ is called the extended Shi arrangement, 
which are known to be free \cite{yos-char}. 
\end{remark}

\begin{theorem} 
\label{thm:main}
Let $\scA$ be a locally $A_2$ arrangement with a positive 
system in $V=\bC^\ell$. We fix a positive system 
$\Phi^+=\{\alpha_H\mid H\in\scA\}\subset V^*$ of defining equations. 
Let $m:\scA\rightarrow\bZ_{\geq 0}$ be a multiplicity. 
We assume the following condition: 
\begin{itemize}
\item[(*)] Let $\scA_X=\{H_i, H_j, H_k\}$ be a 
codimension two localization with 
$\alpha_i=\alpha_j+\alpha_k$. If $m(H_i)$ is 
odd, then at least one of $m(H_j)$ or $m(H_k)$ 
is odd. 
\end{itemize}
Then $(\scA, m)$ is free, if and only if 
it is extendable. 
Indeed, $E(\scA, m)$ is a free arrangement. 
\end{theorem}
We will give the proof in the next section. 
Here we notice an immediate corollary. 

\begin{corollary}
Let $\scA$ be a locally $A_2$ arrangement with a positive 
system. Suppose that the multiplicity $m$ satisfies 
either 
$m(H)$ is odd $\forall H\in\scA$ or 
$m(H)$ is even $\forall H\in\scA$. 
If the multiarrangement $(\scA, m)$ is free, then 
it is extendable. 
\end{corollary}

\begin{remark}
The condition (*) in Theorem \ref{thm:main} is related 
to the following phenomenon. Consider a multiarrangement 
$x^2y^2(x+y)^1$. Then (deconing of) our extension 
$\mathbf{d}E(x^2y^2(x+y)^1)$ 
is defined by 
$$
x(x-1)y(y-1)(x+y), 
$$
which is not free. However another extension 
$$
x(x-1)y(y-1)(x+y-1)
$$
is free. This shows that even $E(\scA, m)$ is not 
free, $(\scA, m)$ might have another free extension. 
The author does not know whether if the extendability 
can be proved without assuming condition (*). 
See for a little more complicated example. 
\end{remark}

\begin{example}
Let us consider a multiarrangement 
$x^4y^4z^4(x+y)^5(y+z)^5(x+y+z)^4$. It is known to 
be free with exponents $(8,9,9)$ (see 
\cite{ay-quasi, yos-prim} or 
Proposition \ref{prop:ay} below). The extension 
$E(x^4y^4z^4(x+y)^5(y+z)^5(x+y+z)^4)$ is defined by 
\begin{eqnarray*}
x,y,z    &=&kw\ (k=     -1, 0, 1, 2)\\
x+y, y+z &=&kw\ (k= -2, -1, 0, 1, 2)\\
x+y+z    &=&kw\ (k=     -1, 0, 1, 2)\\
        w&=&0, 
\end{eqnarray*}
which is not free (look at the localization at 
$x=y=w=0$ and use Lemma \ref{lem:rk2} (3-ii)). 
However another extension 
\begin{eqnarray*}
x,y,z    &=&kw\ (k=     -1, 0, 1, 2   )\\
x+y, y+z &=&kw\ (k=     -1, 0, 1, 2, 3)\\
x+y+z    &=&kw\ (k=         0, 1, 2, 3)\\
        w&=&0, 
\end{eqnarray*}
is free. 
\end{example}

We can check the following for $\ell=3$. 
\begin{question}
Suppose $\scA$ is of type $A_\ell$ and $(\scA, m)$ is 
free. Then is $(\scA, m)$ always extendable? 
\end{question}

\section{Proof} 

Proof of Theorem \ref{thm:main} is done by 
the induction on the rank $\ell$. If $\ell=2$, then 
$\scA$ is either $|\scA|=2$ or type $A_2$. 
Suppose $|\scA|=2$. Then $E(\scA, m)$ is obviously 
free. Suppose $(\scA, m)$ is defined by 
$x^ay^b(x+y)^c$. 
The next lemma is elementary. 
\begin{lemma}\label{lem:rk2}
Assume $a\leq b$. Set $k=a+b+c$ and 
$\scE=E(x^ay^b(x+y)^c)$. 
\begin{itemize}
\item[(1)] If $c<b-a+1$, then 
$\chi(\scE, t)=(t-1)(t-b)(t-a-c)$. 
\item[(2)] If $c\geq a+b+1$, then 
$\chi(\scE, t)=(t-1)(t-a-b)(t-c)$. 
\item[(3)] $b-a\leq c-1< a+b$, 
\begin{itemize}
\item[(i)] $(a, b, c)\neq (\mbox{even, even, odd}) $, then 
$\chi(\scE, t)=(t-1)(t-\lfloor k/2\rfloor )(t-\lceil k/2\rceil )$. 
\item[(ii)] $(a, b, c)= (\mbox{even, even, odd}) $, 
then 
$\chi(\scE, t)=
(t-1)\left((t-\frac{k}{2})^2+\frac{3}{4}\right)$. 
\end{itemize}
\end{itemize}
\end{lemma}

The next result is due to Wakamiko. 

\begin{proposition}\cite{waka-exp} 
Let $(\scA, m)=x^ay^b(x+y)^c$. 
Assume $a\leq b$ and set $k=a+b+c$ as above. 
Since it is rank two, $(\scA, m)$ is always free. 
The exponents are given as follows: 
\begin{itemize}
\item[(1)] If $c<b-a+1$, then 
$\exp(\scA, m)=(b,a+c)$. 
\item[(2)] If $c\geq a+b+1$, then 
$\exp(\scA, m)=(c,a+b)$. 
\item[(3)] $b-a\leq c-1< a+b$, 
then 
$\exp(\scA, m)=(\lfloor k/2\rfloor, \lceil k/2\rceil )$. 
\end{itemize}
\end{proposition}
In \cite{yos-3arr}, a characterization of freeness for 
rank three arrangemenets is given. It can be stated 
as follows. 

\begin{proposition}
For $\ell=2$, 
$E(\scA, m)$ is free with exponents $(1, d_1, d_2)$ if and 
only if 
\begin{itemize}
\item $\chi(E(\scA, m), t)=(t-1)(t-d_1)(t-d_2)$ and 
\item $\exp (\scA, m)=(d_1, d_2)$. 
\end{itemize}
\end{proposition}
Combining these results, we can prove Theorem \ref{thm:main} 
for $\ell=2$. (Note that the condition (*) in the theorem 
is corresponding to that Lemma \ref{lem:rk2} (3) (ii) does 
not occur. ) 

We now consider the case $\ell\geq 3$. 
Let us first recall the following result. 
\begin{proposition}
\label{prop:free}
\cite{yos-char} 
$E(\scA, m)$ is free with exponents $(1, e_1, \ldots, e_\ell)$ 
if and only if $(\scA, m)$ 
is free with exponents $(e_1, \ldots, e_\ell)$ 
and $E(\scA, m)_X$ is free for any positive 
dimensional intersection $X\subset V$. 
\end{proposition}
It is easily seen that 
$
E(\scA, m)_X=
E(\scA_X, m|_{\scA_X})$. 
Since the localization $(\scA_X, m|_{\scA_X})$ of a 
free multiarrangement $(\scA, m)$ is free with rank 
at most $\ell-1$, 
it follows from the inductive 
hypothesis that $E(\scA, m)_X$ is free. 
Hence Propositon \ref{prop:free} shows that 
$E(\scA, m)$ is free. 
\qed


\section{Totally non-free arrangements}

In a recent paper \cite{atw-emulti} Abe, Terao and Wakefield observed 
several phenomena concerning multiplicities and freeness of 
a multiarrangement $(\scA, m)$. In particular they prove that 
generic four planes $(\scA, m)$ defined by 
$x_1^{m_1}x_2^{m_2}x_3^{m_3}(x_1+x_2+x_3)^{m_4}$ 
will never free for any positive multiplicity 
$m:\scA\rightarrow\bZ_{>0}$. 
Such an arrangement $\scA$ is called totally non-free. 
As an application of extendability techniques, 
we give a straightforward proof of 
totally non-freeness for generic arrangements. 

\begin{proposition}
\label{prop:nonfree}
Suppose $\ell=\dim V\geq 3$ and 
$\scA$ is a generic arrangement with $|\scA|>\ell$. 
Let $m:\scA\rightarrow\bZ_{>0}$. 
Then $(\scA, m)$ is not free. 
\end{proposition}

\proof 
Fix a defining equations $\alpha_H$ for each $H$. 
As is already noticed, it is a positive system. 
Since $(\scA, m)$ is locally Boolean, 
$E(\scA, m)_X=E(\scA_X, m|_{\scA_X})$ is free for any 
nonzero subspace $X\subset V\times\{0\}$. Hence if $(\scA, m)$ is 
free, then Proposition \ref{prop:free} shows that 
$E(\scA, m)$ is also free. 
However, let us consider the restriction to 
the subspace $X=\{0\}\times\bC\subset V\times\bC$. Then 
the localization $E(\scA, m)_X$ is isomorphic to $\scA$ which is 
not free. This is a contradiction. \qed

\section{Free interpolations between 
extended Shi and Catalan arrangements}

Let $\scA$ be a crystallographic Coxeter arrangement with 
a fixed positive system $\Phi^+$ of roots. As is already 
mentioned, $E(\scA, 2k+1)$ and $E(\scA, 2k)$ are 
free for any $k\in\bZ_{>0}$. Obviously 
these two families of arrangements are 
related to each other as 
$$
\cdots\subset E(\scA, 2k-1) 
\subset E(\scA, 2k) \subset 
E(\scA, 2k+1)\subset\cdots. 
$$
In \cite{yos-preprint}, it is observed that 
there exist many free arrangements 
$\scB$ such that $E(\scA, 2k)\subset\scB\subset 
E(\scA, 2k+1)$. The purpose of this section is 
to give a complete description of free 
arrangements interpolating these families for type ADE. 

Let $m:\scA\rightarrow\{0,1\}$ be a $\{0,1\}$-valued 
multiplicity. Any interpolating arrangement can 
be described as $E(\scA, 2k\pm m)$ for some $m$. 
We will describe free interpolations by using 
$\{0,1\}$-valued 
multiplicity $m$. 
Our main result in this section is 
the following. 

\begin{theorem}\label{thm:interp}
Let $\scA$ be an irreducible Coxeter arrangement of type ADE with 
the Coxeter number $h$. 
Fix $\Phi^+$ a positive root system. Let $k$ be a positive 
integer. 
Then the following 
conditions are equivalent. 
\begin{itemize}
\item[(1)] $m:\scA\rightarrow \{0,1\}$ satisfies the following 
condition. 
\begin{itemize}
\item[(1-i)] $m^{-1}(1)\subset\scA$ is a free subarrangement 
with exponents $(e_1, \ldots, e_\ell)$. 
\item[(1-ii)] if $\alpha_1=\alpha_2+\alpha_3 (\alpha_i\in\Phi^+)$ and 
$m(H_1)=1$, then at least $m(H_2)=1$ or $m(H_3)=1$. 
\end{itemize} 
\item[(2)] $E(\scA, 2k+m)$ is free with exponents 
$(1, kh+e_1, \ldots, kh+e_\ell)$. 
\item[(3)] $E(\scA, 2k-m)$ is free with exponents 
$(1, kh-e_1, \ldots, kh-e_\ell)$. 
\end{itemize}
\end{theorem}

Before going proof of Theorem \ref{thm:interp}, 
let us recall a result from \cite{ay-quasi}. 

\begin{proposition}\label{prop:ay}
\cite[Corollary 12]{ay-quasi} 
Let $\scA$ be the Coxeter arrangement with the 
Coxeter number $h$, and 
$m:\scA\rightarrow\{0,1\}$ be a multiplicity. 
Let $k\in\bZ_{> 0}$. 
Then the following conditions are equivalent. 
\begin{itemize}
\item $(\scA, m)$ is free with exponents $(e_1, \ldots, e_\ell)$. 
\item $(\scA, 2k+m)$ is free 
with exponents $(kh+e_1, \ldots, kh+e_\ell)$. 
\item $(\scA, 2k-m)$ is free. 
with exponents $(kh-e_1, \ldots, kh-e_\ell)$. 
\end{itemize}
\end{proposition}

First we prove (1)$\Rightarrow$(2). Suppose $m$ 
satisfies (1-ii). 
Then the multiplicity $2k+m$ 
satisfies the condition (*) in Theorem \ref{thm:main}. 
Thus the extension $E(\scA, 2k+m)$ is free if and only if 
the multiarrangement $(\scA, 2k+m)$ is free. But this is 
done by the assumption (1-i) and 
Proposition \ref{prop:ay}. 

The implication (1)$\Rightarrow$(3) is similar. 

Finally let us prove (2)$\Rightarrow$(1). Suppose 
$E(\scA, 2k+m)$ is free. Then by restricting to $H_0$, 
we have (by Theorem \ref{thm:ziegler}), the multiarrangement 
$(\scA, 2k+m)$ is free. Again from Proposition \ref{prop:ay}, 
we have $(\scA, m)$ is free, in other words, 
$m^{-1}(1)\subset\scA$ is a free subarrangement. 
Thus we have (1-i). To prove (1-ii), suppose that 
there exists $H_1$ such that $\alpha_1=\alpha_2+\alpha_3$ 
and $m(H_1)=1, m(H_2)=m(H_3)=0$. 
Then set $X:=H_1\cap H_2\cap H_3$, which is a codimension 
two subspace. 
From Lemma \ref{lem:rk2} (3-ii), the localization 
$E(\scA, 2k+m)_X$ is not free. It is a contradiction. 
Thus (1-ii) is satisfied. 
\qed

Using Terao's factorization theorem, we obtain 
the following corollary. 

\begin{corollary}
Let $\scA$ be a Coxeter arrangement with the 
Coxeter number $h$ and $m:\scA\rightarrow\{0,1\}$ be 
a multiplicity satisfying the condition (1) of 
Theorem \ref{thm:interp}. Then 
$$
\chi(\mathbf{d}E(\scA, 2k\pm m), t)=
\prod_{i=1}^\ell(t-kh\mp e_i). 
$$
\end{corollary}
The above formula implies 
\begin{equation}
\label{eq:FE}
\chi(\mathbf{d}E(\scA, 2k- m), t)=(-1)^\ell
\chi(\mathbf{d}E(\scA, 2k+ m), 2kh-t). 
\end{equation}
We should note that the formula 
(\ref{eq:FE}) 
is very similar to the 
``functional equation'' discovered by 
Postnikov and Stanley \cite{ps-def}. 
It might be worth asking whether 
the formula (\ref{eq:FE}) holds for any 
crystallographic Coxeter arrangement $\scA$ and 
any multiplicity 
$m:\scA\rightarrow\{0,1\}$. 

\begin{remark}
Recently Abe, Nuida and Numata obtained more 
general results for type $A_\ell$ arrangements 
\cite{ann-bic, nui-char}. 
Their results suggest that 
(\ref{eq:FE}) 
holds even for wider class of multiplicities, namely, 
$m:\scA\rightarrow\{-1, 0,1\}$. 
\end{remark}

\medskip

{\bf Acknowledgement.} 
The author would like to thank 
Takuro Abe and Max Wakefield for useful conversations 
and comments. He also thank to the referee for 
pointing out a crucial mistake in the draft and 
giving many suggestions.

\end{document}